%% file: main.tex
\documentclass[graybox]{svmult}

\usepackage{type1cm}        
\usepackage{makeidx}         
\usepackage{graphicx}        
\usepackage{multicol}        
\usepackage[bottom]{footmisc}
\usepackage{bm}
\usepackage[caption=false]{subfig}

\usepackage{newtxtext}       %
\usepackage{newtxmath}       

\makeindex             

\begin{document}

\title*{A parallel solver for a preconditioned space-time boundary element method for the heat equation}
\titlerunning{A parallel solver for a preconditioned space-time BEM for the heat equation}
\author{Stefan Dohr, Michal Merta, G\"unther Of, Olaf Steinbach and Jan Zapletal}
\authorrunning{S. Dohr, M. Merta, G. Of, O. Steinbach and J. Zapletal}
\institute{Stefan Dohr, Olaf Steinbach, G\"unther Of \at Institute of Applied Mathematics, Graz University of Technology,\\Steyrergasse 30, 8010 Graz, Austria \\ \email{stefan.dohr@tugraz.at, o.steinbach@tugraz.at, of@tugraz.at}
\and Michal Merta\textsuperscript{a,b}, Jan Zapletal\textsuperscript{a,b} \at a)~IT4Innovations, V\v{S}B -- Technical University of Ostrava,\\17. listopadu 2172/15, 708 33 Ostrava-Poruba, Czech Republic\\b)~Department of Applied Mathematics, V\v{S}B -- Technical University of Ostrava,\\17. listopadu 2172/15, 708 33 Ostrava-Poruba, Czech Republic\\ \email{michal.merta@vsb.cz, jan.zapletal@vsb.cz}}

\maketitle

\newcommand{\norm}[1]{\left\lVert#1\right\rVert}
\newcommand{\abs}[1]{\ensuremath{\left\vert#1\right\vert}}
\newcommand{\traceH}{H^{1/2, 1/4}(\Sigma)}
\newcommand{\dualtraceH}{H^{-1/2, -1/4}(\Sigma)}

\abstract*{We describe a parallel solver for the discretized weakly singular space-time boundary integral equation of the spatially two-dimensional heat equation. The global space-time nature of the system matrices leads to improved parallel scalability in distributed memory systems in contrast to time-stepping methods where the parallelization is usually limited to spatial dimensions. We present a parallelization technique which is based on a decomposition of the input mesh into submeshes and a distribution of the corresponding blocks of the system matrices among processors. To ensure load balancing, the distribution is based on a cylic decomposition of complete graphs. In addition, the solution of the global linear system requires the use of an efficient preconditioner. We present a robust preconditioning strategy which is based on boundary integral operators of opposite order, and extend the introduced parallel solver to the preconditioned system.}

\section{Introduction}
In this note we describe a parallel solver for the discretized weakly singular space-time boundary integral equation of the spatially two-dimensional heat equation. The global space-time nature of the system matrices leads to improved parallel scalability in distributed memory systems in contrast to time-stepping methods where the parallelization is usually limited to spatial dimensions. We present a parallelization technique which is based on a decomposition of the input mesh into submeshes and a distribution of the corresponding blocks of the system matrices among processors. To ensure load balancing, the distribution is based on a cyclic decomposition of complete graphs \cite{KravcenkoMertaZapletal2018, LukasKovarKovarovaMerta2015}. In addition, the solution of the global linear system requires an efficient preconditioner. We present a robust preconditioning strategy which is based on boundary integral operators of opposite order \cite{Hiptmair:2006, SteinbachWendland:1998}.

The parallelization of the discretized space--time integral equation in distributed and shared memory is discussed in \cite{DohrZapMertaOf}. Here, we extend the parallel solver to the preconditioned system. We demonstrate the method for the spatially two-dimensional case. However, the presented results, particularly the parallelization in distributed memory and the stability results for the preconditioner, can be used to extend the method to the three-dimensional problem.

Let $\Omega \subset \mathbb{R}^{2}$ be a bounded domain with a Lipschitz boundary 
$\Gamma := \partial \Omega$ and $T > 0$. 
As a model problem we consider the initial Dirichlet boundary value problem for 
the heat equation
\begin{equation}
\label{eq:model_problem}
\alpha \partial_{t} u - \upDelta_{x}u = 0 \; \textrm{in} \; 
Q:= \Omega \times (0,T), \;
u = g \; \textrm{on} \; \Sigma := \Gamma \times (0,T),\;
u = u_{0} \; \textrm{in} \; \Omega
\end{equation}
with a heat capacity constant $\alpha > 0$, the given initial datum 
$u_{0} \in L^{2}(\Omega)$, and the boundary datum $g \in \traceH$. 
An explicit formula for the solution of \eqref{eq:model_problem} is given by the representation formula for the heat equation \cite{ArnoldNoon:1987}, i.e. for $(x,t) \in Q$ we have
\begin{equation}
\label{eq:representation_formula}
\begin{aligned}
u(x,t)
&=  (\widetilde{M}_0 u_0)(x,t) + (\widetilde{V}w)(x,t) - (Wg)(x,t)\\
&= \int_{\Omega} U^{\star}(x-y, t) u_{0}(y) \,\D y + \frac{1}{\alpha} \int_{\Sigma} U^{\star}(x-y, t-\tau) w(y,\tau) \,\D s_{y} \,\D \tau\\
&\quad - \frac{1}{\alpha} \int_{\Sigma} \frac{\partial}{\partial n_{y}} U^{\star}(x-y, t-\tau) g(y,\tau) \,\D s_{y} \,\D \tau,
\end{aligned}
\end{equation}
with $w := \partial_{n}u$ and $U^{\star}$ denoting the fundamental solution of the two-dimensional heat equation given by
\begin{equation*}
U^{\star}(x-y,t-\tau) =
  \begin{cases}
   \displaystyle \; \frac{\alpha}{4 \pi (t-\tau)} \exp{\left(\frac{-\alpha|x-y|^{2}}{4 (t-\tau)}\right)}\; & \text{for } \tau < t, \\
   \displaystyle \; 0        &\text{otherwise}.
  \end{cases}
\end{equation*}
The yet unknown Neumann datum $w \in \dualtraceH$ can be found by applying the interior Dirichlet trace operator to \eqref{eq:representation_formula} and solving the resulting weakly singular boundary integral equation
\begin{equation}
\label{eq:first_bie}
g(x,t) = (M_{0}u_{0})(x,t) + (Vw)(x,t) + ((\frac{1}{2}I-K)g)(x,t) 
\quad \textrm{for } (x,t) \in \Sigma.
\end{equation}
The operators in \eqref{eq:first_bie} are obtained by composition of the heat potentials in \eqref{eq:representation_formula} with the Dirichlet trace operator.
The ellipticity \cite{costabel:1990} and boundedness of the single-layer operator $V\colon \dualtraceH \rightarrow \traceH$ together with the boundedness of the double-layer operator $K\colon \traceH \rightarrow \traceH$ and the initial Dirichlet operator $M_{0}\colon L^{2}(\Omega) \rightarrow \traceH$ ensure unique solvability of \eqref{eq:first_bie}.

\begin{figure}[htb]
\centering
\subfloat[Tensor product decomposition. \label{fig:sample_boundary_decomposition_2D}]{%
\hspace{0.075\textwidth}
\includegraphics[width=0.35\textwidth]{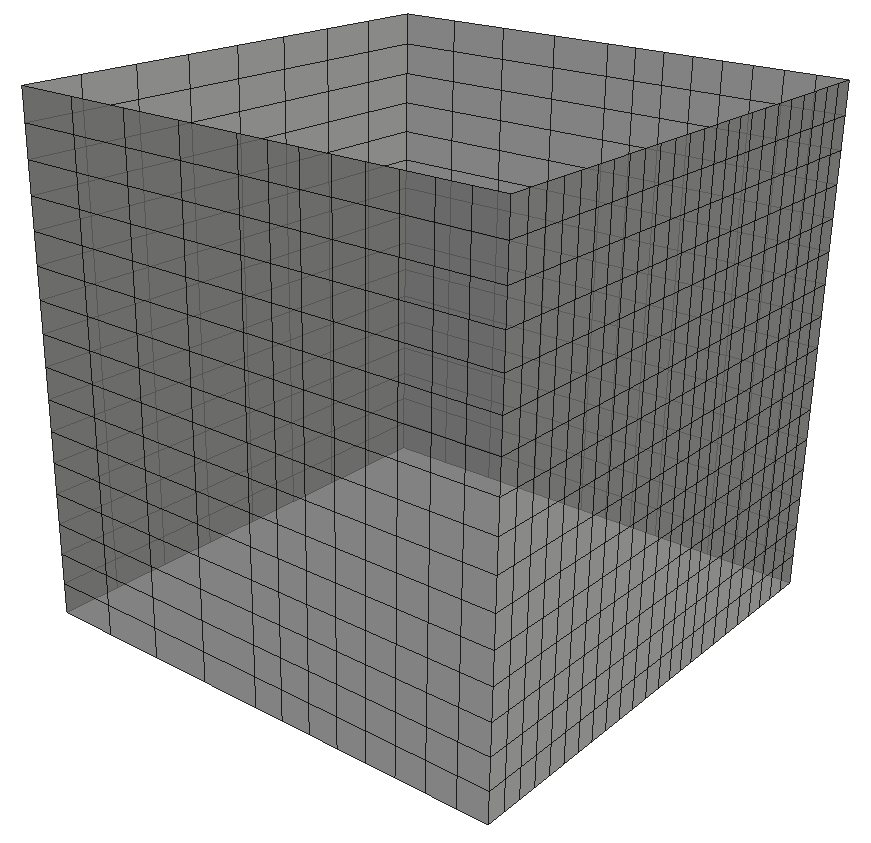}}
\hfill
\subfloat[Submeshes. \label{fig:sample_submeshes}]{%
\includegraphics[width=0.35\textwidth]{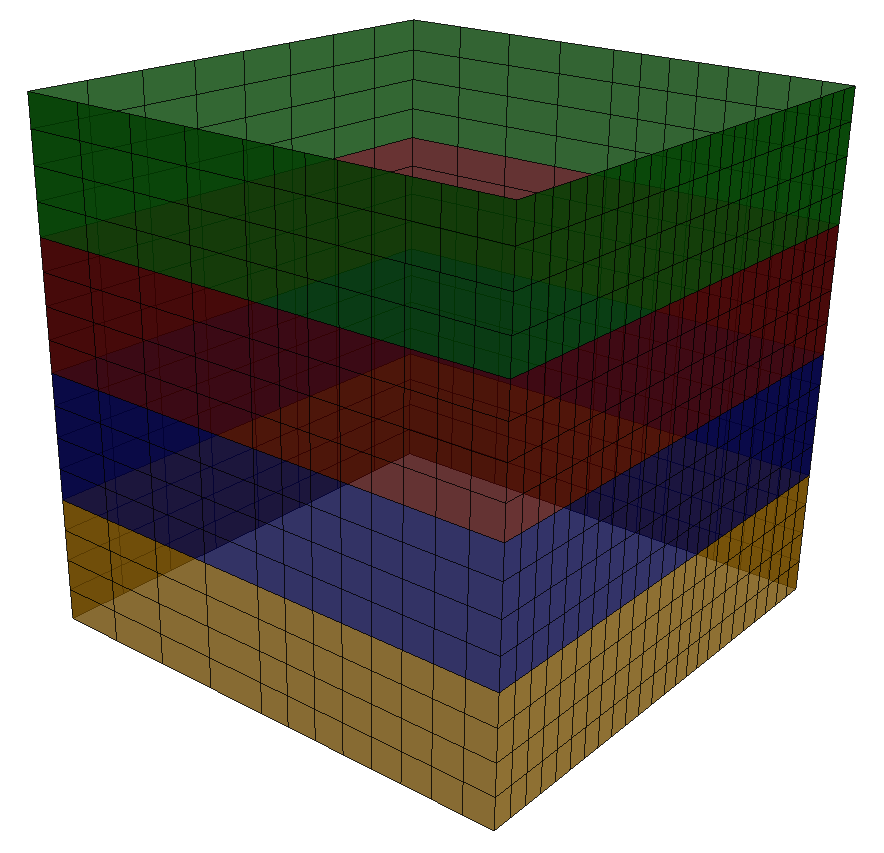}}
\hspace{0.075\textwidth}
\caption{Sample space-time boundary decompositions for $Q = (0,1)^{3}$.}
\label{fig:boundary_decompositions_2D}
\end{figure}

We consider a space-time tensor product decomposition of $\Sigma$~\cite{costabel:1990, MessnerSchanzTausch14, noon:1988} and use the Galerkin method for the discretization of \eqref{eq:first_bie}. For a triangulation $\Gamma_{h} = \left\{\gamma_{i}\right\}_{i=1}^{N_{\Gamma}}$ of the boundary $\Gamma$ and a decomposition $I_{h} = \left\{\tau_{k}\right\}_{k=1}^{N_{I}}$ of the time interval $I := (0,T)$ we define $\Sigma_{h} := \left\{ \sigma = \gamma_{i} \times \tau_{k}\colon i=1,...,N_{\Gamma};\, k=1,...,N_{I}\right\}$, i.e. $\Sigma_h = \left\{\sigma_\ell\right\}_{\ell = 1}^{N}$ with $N = N_\Gamma N_I$. In the two-dimensional case the space-time boundary elements $\sigma$ are rectangular. A sample decomposition of the space-time boundary of $Q = (0,1)^{3}$ is shown in Fig.~\ref{fig:sample_boundary_decomposition_2D}.

We use the space $X_{h}^{0,0}(\Sigma_{h}) := 
\textrm{span}\big\{\varphi_{\ell}^{0}\big\}_{\ell=1}^{N}$
of piecewise constant basis functions and the space $X_{h}^{1,0}(\Sigma_{h}) := \textrm{span} \big\{\varphi_{i}^{10}\big\}_{i=1}^{N}$ of functions that are piecewise linear and globally continuous in space and piecewise constant in time for the approximations of the Cauchy data $w$ and $g$, respectively. The initial datum $u_{0}$ is discretized by using the space of piecewise linear and globally continuous functions $S_{h}^{1}(\Omega_{h}) := \textrm{span}\big\{\varphi_{j}^{1}\big\}_{j=1}^{M_\Omega}$, which is defined with respect to a given triangulation $\Omega_{h} := \big\{\omega_{i}\big\}_{i=1}^{N_{\Omega}}$ of the domain~$\Omega$. This leads to the system of linear equations 
\begin{equation}
\label{eq:linear_equations}
V_{h} \vec{w} = (\frac{1}{2}M_{h} + K_{h})\vec{g} - M_{h}^{0}\vec{u_{0}}
\end{equation}
 with
\begin{equation*}
\begin{aligned}
V_{h}[\ell, k] &:= \langle V \varphi_k^{0}, \varphi_{\ell}^{0} \rangle_{L^{2}(\Sigma)}, & K_{h}[\ell, i] &:= \langle K \varphi_i^{10}, \varphi_{\ell}^{0} \rangle_{L^{2}(\Sigma)}, \\
M_{h}^{0}[\ell, j] &:=\langle M_0 \varphi_j^{1}, \varphi_{\ell}^{0} \rangle_{L^{2}(\Sigma)}, & M_{h}[\ell, i] &:= \langle \varphi_i^{10}, \varphi_{\ell}^{0} \rangle_{L^{2}(\Sigma)},
\end{aligned}
\end{equation*}
for $i, k, \ell = 1,\ldots,N$ and $j=1,\ldots,M_\Omega$. Due to the ellipticity of the single-layer operator $V$ the matrix 
$V_{h}$ is positive definite and therefore \eqref{eq:linear_equations} is uniquely solvable.

\section{Operator preconditioning}
The boundary element discretization is done with respect to the
whole space-time boundary $\Sigma$ and since we want to solve \eqref{eq:linear_equations} without an application of time-stepping schemes to make use of parallelization in time, we need to develop an efficient
iterative solution technique. The linear system \eqref{eq:linear_equations} with the positive definite but non-symmetric matrix $V_h$ can be solved by a
preconditioned GMRES method. Here we will apply a preconditioning technique 
based on boundary integral operators of opposite order 
\cite{SteinbachWendland:1998}, also known as operator or 
Calderon preconditioning \cite{Hiptmair:2006}. 

First, we introduce the hypersingular operator $D$, which is defined as the negative Neumann trace of the double layer potential $W$ in \eqref{eq:representation_formula}, i.e. $(Dv)(x,t) = - \partial_n (Wv)(x,t)$ for $(x,t) \in \Sigma$.
The single-layer
operator $V\colon \dualtraceH \to \traceH$ and the hypersingular operator $D\colon \traceH \to \dualtraceH$ are both elliptic \cite{costabel:1990} and the composition $DV\colon \dualtraceH \to \dualtraceH$ defines an operator of order zero. 
Thus, following \cite{Hiptmair:2006}, the Galerkin discretization of $D$ allows the construction of a suitable preconditioner for $V_h$.
While the discretization of the single-layer operator $V$ is done with respect to $X_h^{0,0}(\Sigma_h)$, for the Galerkin discretization of the hypersingular
operator $D$ we need to use a conforming trial space
$Y_{h} = \textrm{span} \left\{\psi_i\right\}_{i=1}^{N} \subset \traceH$, see also \cite{DohrSteinbach2018} for the spatially one-dimensional problem.

\begin{theorem}[\cite{Hiptmair:2006,SteinbachWendland:1998}]
Assume that the discrete stability condition
\begin{equation}
\label{eq:stability_trial_spaces}
\sup_{0 \neq v_{h} \in Y_{h}} \frac{\langle \tau_{h}, v_{h}
  \rangle_{L^{2}(\Sigma)}}{\norm{v_{h}}_{\traceH}} \geq c_{1}^{M}
\norm{\tau_{h}}_{\dualtraceH} \quad \textrm{for all } \tau_{h} \in X_h^{0,0}(\Sigma_h)
\end{equation}
holds. Then there exists a constant $c_{\kappa} > 1$ such that $\kappa\left(M_{h}^{-1} D_{h} M_{h}^{-\top} V_{h}\right) \leq c_{\kappa}$
where, for $k,\ell=1,\ldots,N$,
\begin{equation*}
D_{h}[\ell, k] = \langle D \psi_{k}, \psi_{\ell}\rangle_{\Sigma} \;, 
\quad M_{h}[\ell, k] = \langle \varphi_{k}^0, \psi_{\ell} \rangle_{L^{2}(\Sigma)}
\; .
\end{equation*}
\end{theorem}
Thus we can use $C_{V}^{-1} = M_{h}^{-1} D_{h} M_{h}^{-\top}$ as a
preconditioner for $V_h$. For the computation of the matrix $D_h$ we use an alternative representation of the associated bilinear form which is attained by applying integration by parts, see \cite[Theorem 6.1]{costabel:1990}. It remains to define a suitable boundary element space $Y_h$ such that the
stability condition \eqref{eq:stability_trial_spaces} is satisfied.
In what follows we will discuss a possible choice.

We assume that the decompositions $\Gamma_{h}$ and $I_h$ are locally quasi-uniform. For the given boundary element mesh $\Gamma_h$ we construct a dual mesh $\widetilde{\Gamma}_h := \left\{\tilde{\gamma}_{\ell}\right\}_{\ell}^{N_\Gamma}$ according to \cite{HiptJerezUrzua:2014, Steinbach:2002} and assume, that $\widetilde{\Gamma}_h$ is locally quasi-uniform as well. For the discretization of the operator $D$ we choose $Y_h = X_h^{1,0}(\widetilde{\Sigma}_h) \subset \traceH$, which denotes the space of functions that are piecewise linear and globally continuous in space and piecewise constant in time, defined with respect to the decompositions $\widetilde{\Gamma}_h$ and $I_h$, respectively. In order to prove the stability condition \eqref{eq:stability_trial_spaces} we establish the $\traceH$-stability of the $L^{2}(\Sigma)$-projection $\widetilde{Q}_h^{1,0}\colon L^{2}(\Sigma) \rightarrow Y_h \subset L^{2}(\Sigma)$ defined by
\begin{equation}
\label{eq:def_Qh}
\langle \widetilde{Q}_h^{1,0}v, \tau_h \rangle_{L^{2}(\Sigma)} = \langle v, \tau_h  \rangle_{L^{2}(\Sigma)} \quad \textrm{for all } \tau_h \in X_h^{0,0}(\Sigma).
\end{equation}
The Galerkin-Petrov variational problem \eqref{eq:def_Qh} is uniquely solvable since the trial and test spaces satisfy a related stability condition \cite{DohrNiinoSteinbach}. When assuming appropriate local mesh conditions of $\Gamma_h$ and $\widetilde{\Gamma}_h$, see \cite{Steinbach:2001, Steinbach:2002}, we are able to establish the stability of $\widetilde{Q}_h^{1,0}\colon \traceH \rightarrow \traceH$, see \cite{DohrNiinoSteinbach} for a detailed discussion. Hence, there exists a constant $c_{S} > 0$ such that
\begin{equation}
\label{eq:stability_l2_projection}
\norm{\widetilde{Q}_h^{1,0}v}_{\traceH} \leq c_{S} \norm{v}_{\traceH} \quad \textrm{for all } v \in \traceH.
\end{equation}
The stability estimate \eqref{eq:stability_l2_projection} immediately implies the stability condition \eqref{eq:stability_trial_spaces}. Hence the condition number $\kappa(C_V^{-1}V_h)$ with $C_{V}^{-1} = M_{h}^{-1} D_{h} M_{h}^{-\top}$ is bounded.

\section{Distributed memory parallelization}

Distributed memory parallelization of the solver is based on the scheme presented in \cite{KravcenkoMertaZapletal2018, LukasKovarKovarovaMerta2015} for spatial problems. In \cite{DohrZapMertaOf} we have extended the approach to support time-dependent problems for the heat equation. Let us briefly describe the method and refer the more interested readers to the above-mentioned papers.

\begin{figure}[hbt]
\centering
	\subfloat[Generating graph.]{
    \includegraphics[width=0.23\textwidth]{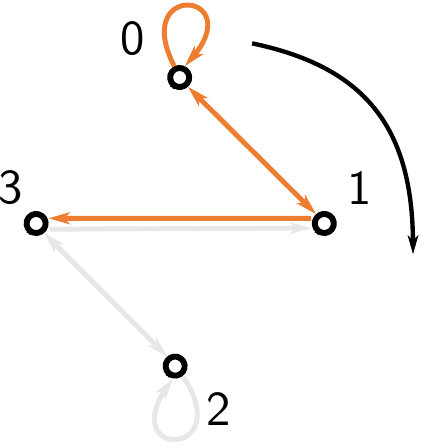}
	}%
     \hspace{0.2\textwidth}
  \subfloat[Block distribution.]{
    \includegraphics[width=0.23\textwidth]{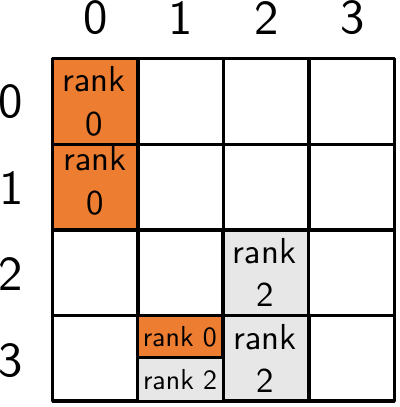}
	}%
	\caption{Distribution of the system matrix blocks among four processes.}
\label{fig:example_decomposition}
\end{figure}

To distribute the system among $P$ processes the space-time mesh $\Sigma$ is decomposed into $P$ slices in the temporal dimension (see Fig.~\ref{fig:sample_submeshes}) which splits the matrices $A\in\{V_h, K_h, D_h\}$ into $P\times P$ blocks 

\begin{equation*}
A = \begin{bmatrix}
A_{0, 0} & 0 & \cdots & 0 \\
A_{1, 0} & A_{1, 1} & \cdots & 0 \\
\vdots & \vdots & \ddots & \vdots \\
A_{P-1, 0} & A_{P-1, 1} & \cdots & A_{P-1, P-1} \\
\end{bmatrix}.
\end{equation*}
The matrices are block lower triangular with lower triangular blocks on the main diagonal due to the properties of the fundamental solution and the selected discrete spaces. We aim to distribute the blocks among processes such that the number of shared mesh parts is minimal and each process owns a single diagonal block. For this purpose we consider each block $A_{i,j}$ as an edge $(i,j)$ of a complete graph $K_P$ on $P$ vertices. The distribution problem corresponds to finding a suitable decomposition of $K_P$ into $P$ subgraphs $G_0, G_1, \ldots, G_{P-1}$. In \cite{DohrZapMertaOf, KravcenkoMertaZapletal2018} we employ a cyclic decomposition algorithm~-- first, a generator graph $G_0$ on a minimal number of vertices (corresponding to blocks to be assembled by the process 0) is constructed; the remaining graphs $G_1, \ldots, G_{P-1}$ are obtained by a clock-wise rotation of $G_0$ along vertices of $K_P$ placed on a circle. An example of the generating graph and the corresponding matrix decomposition for four processes is depicted in Fig.~\ref{fig:example_decomposition}. In the case of the initial matrix $M_h^0$ we distribute block-rows of the matrix among processes. Similarly, since the matrix $M_h$ is block diagonal, each process owns exactly one block of the matrix.

In addition to the distributed memory parallelization by MPI, the assembly of the matrices is parallelized and vectorized in shared memory using OpenMP \cite{DohrZapMertaOf,ZapMerMal2016}. Therefore, in our numerical experiments we usually employ hybrid parallelization using one MPI process per CPU socket and an appropriate number of OpenMP threads per process.

\section{Numerical experiments}

The presented examples refer to the initial Dirichlet boundary value problem \eqref{eq:model_problem} on the space-time domain $Q := (0,1)^{3}$. The heat capacity constant is set to $\alpha = 10$. We consider the exact solution
\begin{equation*}
u(x,t) := \exp{\left(-\frac{t}{\alpha}\right)} \sin{\left(x_1 \cos{\frac{\pi}{8}} + x_2 \sin{\frac{\pi}{8}}\right)} \quad \textrm{for } (x,t) = (x_1,x_2,t) \in Q
\end{equation*}
and determine the Dirichlet datum $g$ and the initial datum $u_0$ accordingly. The linear system~\eqref{eq:linear_equations} is solved by the GMRES method with a relative precision of~$10^{-8}$. 

\paragraph{Operator preconditioning}
As a preconditioner we use the discretization $C_V^{-1} = M_h^{-1}D_h M_h^{-\top}$ of the hypersingular operator $D$ in the space $X_h^{1,0}(\widetilde{\Sigma}_h)$, while the Galerkin discretization of the integral equation \eqref{eq:first_bie} is done with respect to $X_h^{0,0}(\Sigma_h)$. Instead of using $M_h$ in the preconditioner we computed a lumped mass matrix. Thus, the matrix becomes diagonal and the inverse can be applied efficiently.

The example corresponds to a globally uniform boundary element mesh with the mesh size $h = \mathcal{O}(2^{-L})$. Table \ref{tab_prec_uniform} shows the iteration numbers of the non-preconditioned and preconditioned GMRES-method. As expected, the iteration numbers of the preconditioned version are bounded due to the boundedness of $\kappa(C_V^{-1}V_h)$. For numerical results in the case of an adaptive refinement we refer to \cite{DohrNiinoSteinbach}.
\begin{table}[b]
\centering
\setlength{\tabcolsep}{10pt}
\caption{Iteration numbers in the case of 
uniform refinement.}
\label{tab_prec_uniform} 
\input{gmres_iteration_numbers.tex}
\end{table}

\paragraph{Scalability in distributed memory}
The numerical experiments for the scalability were executed on the Salomon cluster at IT4Innovations National Supercomputing Center in Ostrava, Czech Republic. The cluster is equipped with 1008 nodes with two 12-core Intel Xeon E5-2680v3 Haswell processors and 128 GB of RAM. Nodes of the cluster are interconnected by the InfiniBand 7D enhanced hypercube network.

We tested the assembly of the BEM matrix $D_h$. Computation times for the assembly of the matrices $V_h$, $K_h$, $M_h^0$, the related matrix-vector multiplication, and the evaluation of the solution in $Q$ can be found in \cite{DohrZapMertaOf}. Strong scaling of the parallel solver was tested using a tensor product decomposition of $\Sigma$ into $65\,536$, $262\,144$ and $1\,048\,576$ space-time surface elements.
We used up to 256 nodes ($6\, 144$ cores) of the Salomon cluster for the computations and executed two MPI processes per node. Each MPI process used 12 cores for the assembly of the matrix blocks.

In Table \ref{tab:assembly_Dh} the assembly times for $D_h$ including the speedup and efficiency are listed. We obtain almost optimal parallel scalability. Note that the number of nodes is restricted by the number of elements of the temporal decomposition $I_h$. Conversely, for fine meshes we need a certain number of nodes to store the matrices.

\begin{table}[t]
\centering
\caption{Assembly of $D_{h}$ on 65\,536, 262\,144, and 1\,048\,576 space-time elements.}
\label{tab:assembly_Dh}
\input{mpi_Dh.tex}
\end{table}

\section{Conclusion}

In this note we have described a parallel space-time boundary element solver for the two-dimensional heat equation. The solver is parallelized using MPI in the distributed memory. The distribution of the system matrices is based on \cite{DohrZapMertaOf,KravcenkoMertaZapletal2018, LukasKovarKovarovaMerta2015}. The space-time boundary mesh is decomposed into time slices which define blocks in the system matrices. These blocks are distributed among MPI processes using the graph decomposition based scheme. For a detailed discussion on shared memory parallelization see \cite{DohrZapMertaOf}. Moreover, we have introduced an efficient preconditioning strategy for the space-time system which is based on the use of boundary integral operators of opposite order. The preconditioner was then distributed with the presented parallelization technique. 

The numerical experiments for the proposed preconditioing strategy confirm the theoretical findings, i.e. the boundedness of the iteration number of the iterative solver. We also tested the efficiency of the parallelization scheme for the preconditioner. The results show almost optimal scalability.

\acknowledgement{The research was supported by the project `Efficient parallel implementation of boundary element methods' provided jointly by the Ministry of Education, Youth and Sports (7AMB17AT028) and OeAD (CZ~16/2017). SD acknowledges the support provided by the International Research Training Group 1754, funded by the German Research Foundation (DFG) and the Austrian Science Fund (FWF). JZ an MM further acknowledge the support provided by The Ministry of Education, Youth and Sports from the National Programme of Sustainability (NPS II) project `IT4Innovations excellence in science -- LQ1602' and the Large Infrastructures for Research, Experimental Development and Innovations project `IT4Innovations National Supercomputing Center -- LM2015070'.}

\bibliographystyle{plain}
\bibliography{references}

\end{document}

%% file: gmres_iteration_numbers.tex
\begin{tabular}
{rrrr}
\hline\noalign{\smallskip}
$L$ & $N$ & It. & It. prec.  \\
\noalign{\smallskip}\svhline\noalign{\smallskip}
2&64&14&17\\
3&256&19&18\\
4&1\,024&24&20\\
5&4\,096&35&20\\
6&16\,384&50&20\\
7&65\,536&67&20\\
8&262\,144&91&19\\
9&1\,048\,576&122&19\\
\noalign{\smallskip}\hline\noalign{\smallskip}
\end{tabular}

%% file: mpi_Dh.tex
\begin{tabular}{p{1.3cm}p{0.9cm}p{0.9cm}p{0.9cm}p{0.9cm}p{0.9cm}p{0.9cm}p{0.9cm}p{0.9cm}p{0.9cm}}
\hline\noalign{\smallskip}
  nodes $\downarrow$ & \multicolumn{3}{c}{$D_{h}$ assembly [s]} & \multicolumn{3}{c}{$D_{h}$ speedup} & \multicolumn{3}{c}{$D_{h}$ efficiency [\%]} \\
  mesh $\rightarrow$ & 65k & 262k & 1M & 65k & 262k & 1M & 65k & 262k & 1M \\
\noalign{\smallskip}\svhline\noalign{\smallskip}
1 & 184.1 & --- & --- & 1.0 & --- & --- & 100.0 & --- & --- \\
2 & 92.0 & --- & --- & 2.0 & --- & --- & 100.1 & --- &--- \\
4 & 46.8 & --- & --- & 3.9 & --- & --- & 98.4 & --- & --- \\
8 & 23.8 & 373.6 & --- & 7.7 & 1.0 & --- & 96.7 & 100.0 & --- \\
16 & 11.8 & 186.1 & --- & 15.6 & 2.0 & --- & 97.3 & 100.4 & --- \\
32 & 5.9 & 91.9 & --- & 31.0 & 4.1 & --- & 96.7 & 101.7 & --- \\
64 & 3.0 & 47.0 & 747.0 & 60.5 & 7.9 & 1.0 & 94.6 & 99.3 & 100.0 \\
128 & --- & 24.0 & 376.9 & --- & 15.6 &  2.0 & --- & 97.3 & 99.1 \\
256 & --- & --- & 193.5 & --- & --- & 3.9 & --- & --- & 96.5\\
\noalign{\smallskip}\hline\noalign{\smallskip}
\end{tabular}